\def\UseSection{
      \numberwithin{equation}{section}
	\theoremstyle{plain}
      \newtheorem{theorem}    {Theorem}[section]
      \DefineTheorems 
}

\def\DefineTheorems{
	
	\newtheorem{lemma}      [theorem] {Lemma}
	
	\newtheorem{prop}       [theorem] {Proposition}
	
	\newtheorem{cor}        [theorem] {Corollary}

	\theoremstyle{definition}
	\newtheorem{defn}       [theorem] {Definition}

	\theoremstyle{definition}

}

\newcommand{\bt}   {\begin{theorem}}
\newcommand{\et}   {\end  {theorem}}
\newcommand{\bl}   {\begin{lemma}}
\newcommand{\el}   {\end  {lemma}}
\newcommand{\bp}   {\begin{prop}}
\newcommand{\ep}   {\end  {prop}}
\newcommand{\bc}   {\begin{cor}}
\newcommand{\ec}   {\end  {cor}}
\newcommand{\bd}   {\begin{defn}}
\newcommand{\ed}   {\end  {defn}}

\newcommand{\ba}   {\begin{array}}
\newcommand{\ea}   {\end  {array}}
\newcommand{\be}   {\begin{enumerate}}
\newcommand{\ee}   {\end  {enumerate}}
\newcommand{\bi}   {\begin{itemize}}
\newcommand{\ei}   {\end  {itemize}}

\def\eq#1\en{\begin{equation}#1\end{equation}}  
\def\eqsplit#1\ensplit{
	\begin{equation}\begin{split}#1\end{split}\end{equation}
	}
\def\eqalign#1\enalign{
	\begin{align}#1\end{align}
	}
\def\eqmul#1\enmul{
	\begin{multline}#1\end{multline}
	}
\newcommand{\eqarrstar} {\begin{eqnarray*}} 
\newcommand{\enarrstar} {\end{eqnarray*}} 
\newcommand{\eqarray}   {\begin{eqnarray}} 
\newcommand{\enarray}   {\end{eqnarray}} 
 
\newcommand{\smallE}{\scriptstyle \rightarrow}
\newcommand{\smallW}{\scriptstyle \leftarrow}

\newcommand{\smallS}{\scriptstyle \downarrow}

\makeatletter
\newcommand{\labelcounter}[2]{{%
	\stepcounter{#1}
	\protected@write\@auxout{}%
	{\string\newlabel{#2}{{\csname the#1\endcsname}{\thepage}}}%
	{\ref{#2}}
	}}
\makeatother
\newcommand{\sss}   { \scriptscriptstyle } 

\newcommand{\Zbold} {{\mathbb Z}}

\newcommand{\spose}[1] {{\hbox to 0pt{#1\hss}} }
\newcommand{\ltapprox} {\mathrel{\spose{\lower 3pt\hbox{$\mathchar"218$}}
\raise 2.0pt\hbox{$\mathchar"13C$}}}
\newcommand{\gtapprox} {\mathrel{\spose{\lower 3pt\hbox{$\mathchar"218$}}
\raise 2.0pt\hbox{$\mathchar"13E$}}}

\documentclass[12pt]{amsart}
\usepackage{amsmath,amssymb,amsfonts,amsthm,bbm}
\usepackage[dvips]{graphicx}
\newtheorem{THM}{Theorem}[section]

\newtheorem{LEM}[THM]{Lemma}

\newtheorem{DEF}[THM]{Definition}

\UseSection  
\newcommand{\hlf}{\frac{1}{2}}
\newcommand{\ra}{\rightarrow}

\renewcommand{\to}      {\rightarrow}

\newcounter{countC}  
\newcounter{countR}  
\newcommand{\Z}{\Zbold}

\newcommand{\mc}[1]{\mathcal{#1}}
\newcommand{\mP}{\mathbb{P}}

\newcommand{\mE}{\mathbb{E}}

\newcommand{\LR}{\leftrightarrow}
\newcommand{\WE}{\LR}

\newcommand{\NE}{\begin{picture}(,)
\put(2,-5){$\rightarrow$}
\put(0,.5){$\uparrow$}
\end{picture}\hspace{.5cm}
}

\newcommand{\NEalpha}{\begin{picture}(,)
\put(2,-5){$\rightarrow$}
\put(4,-8){$\scriptstyle\alpha$}
\put(0,.5){$\uparrow$}
\end{picture}\hspace{.5cm}
}
\newcommand{\NEbeta}{\begin{picture}(,)
\put(2,-5){$\rightarrow$}
\put(4,-9){$\scriptstyle\beta$}
\put(0,.5){$\uparrow$}
\end{picture}\hspace{.5cm}
}
\newcommand{\SE}{\begin{picture}(,)
\put(1.5,4.8){$\rightarrow$}
\put(-.5,-.5){$\downarrow$}
\end{picture}\hspace{.5cm}
}

\newcommand{\SW}{\begin{picture}(,)
\put(0,4.8){$\leftarrow$}
\put(8.2,-0.5){$\downarrow$}
\end{picture}\hspace{.5cm}
}

\newcommand{\NW}{\begin{picture}(,)
\put(0,-5){$\leftarrow$}
\put(8,0){$\uparrow$}
\end{picture}\hspace{.5cm}
}

\newcommand{\NSE}{\begin{picture}(,)
\put(0,0){$\updownarrow$}
\put(2.2,0){$\rightarrow$}
\end{picture}\hspace{.5cm}
}
\newcommand{\SWE}{\begin{picture}(,)
\put(0,5){$\leftarrow$}
\put(5,5){$\rightarrow$}
\put(5.5,-0.5){$\downarrow$}
\end{picture}\hspace{.6cm}
}
\newcommand{\smallSWE}{\begin{picture}(,)
\put(0,2.5){$\smallW$}
\put(3,2.5){$\smallE$}
\put(3.8,-1.5){$\smallS$}
\end{picture}\hspace{.6cm}
}
\newcommand{\NWE}{\begin{picture}(,)
\put(0,-5){$\leftarrow$}
\put(5,-5){$\rightarrow$}
\put(5.5,0.5){$\uparrow$}
\end{picture}\hspace{.6cm}
}

\newcommand{\NSEWalt}{\begin{picture}(,)
\put(8,0){$\updownarrow$}
\put(0,0){$\longleftrightarrow$}
\end{picture}\hspace{.5cm}
}

\newcommand{\blank}[1]{}

\newcommand{\Qed}{\qed \medskip}

\usepackage[usenames]{color}

\begin{document}

\title[RWDRE]{Speed calculations for random walks in degenerate random environments.}

\author[Holmes]{Mark Holmes}
\address{Department of Statistics, University of Auckland}
\email{mholmes@stat.auckland.ac.nz}
\author[Salisbury]{Thomas S. Salisbury} 
\address{Department of Mathematics and Statistics, York University}
\email{salt@yorku.ca}

\maketitle

\begin{abstract}
We calculate explicit speeds for random walks in uniform degenerate random environments. For certain non-uniform random environments, we calculate speeds that are non-monotone. 
\end{abstract}

\section{Introduction}

In \cite{HS_RWDRE} the authors study random walk in an IID random environment, where the environment need not satisfy any ellipticity condition. In other words, where various nearest neighbour transitions may have quenched probability $=0$. If such a walk can get stuck on a finite set of vertices with positive probability, then it will get stuck with probability one. There are necessary and sufficient conditions for such a walk not to get stuck in this way, and \cite{HS_RWDRE} studies transience and speed questions for such walks. There are many interesting models in which such properties are non-trivial. There are also examples in which transience is essentially trivial, but in which speeds can be calculated explicitly, because of a renewal structure. \cite{HS_RWDRE} gives a table of such speeds, for random walks in {\em uniform} degenerate random environments. That is, environments in which the walker chooses at random from the (random) set of allowed steps. The purpose of this note is to supply details for the latter calculations. We also include details of some calculations, for speeds and other quantities, related to examples of non-monotone behaviour.

For fixed $d\ge 2$ let $\mc{E}=\{\pm e_i: i=1,\dots,d\}$ be the set of unit vectors in $\Z^d$.  Let $\mc{P}=M_1(\mc{E})$ denote the set of probability measures on $\mc{E}$, and let $\mu$ be a probability measure on $\mc{P}$.  Let $\Omega=\mc{P}^{\Z^d}$ be equipped with the product measure $\nu=\mu^{\otimes \Z^d}$ (and the corresponding product $\sigma$-algebra).  A random environment $\omega=(\omega_x)_{x\in \Z^d}$ is an element of $\Omega$.  We write $\omega_x(e)$ for $\omega_x(\{e\})$. Note that $(\omega_x)_{x\in \Z^d}$ are i.i.d.~with law $\mu$ under $\nu$.

The random walk in environment $\omega$ is a time-homogeneous Markov chain with transition probabilities from $x$ to $x+e$ defined by 
\begin{equation}
\label{eq:trans_prob}
p_{\omega}(x,x+e)=
\omega_x(e).
\end{equation}
Given an environment $\omega$, we let $\mP_{\omega}$ denote the law of this random walk $X_n$, starting at the origin.  Let $P$ denote the law of the annealed random walk, i.e.~$P(\cdot, \star):=\int_{\star}\mP_{\omega}(\cdot)d\nu$.
Since $P(A)=E_{\nu}[\mP_{\omega}(A)]$ and $0\le f(\omega)=\mP_{\omega}(A)\le 1$, $P(A)=1$ if and only if $\mP_{\omega}(A)=1$ for $\nu$-almost every $\omega$.  Similarly $P(A)=0$ if and only if $\mP_{\omega}(A)=0$ for $\nu$-almost every $\omega$.  
If we start the RWRE at $x\in\Z^d$ instead, we write $P_x$ for the corresponding probability, so $P=P_o$.

We associate to each environment $\omega$ a directed graph $\mc{G}(\omega)$ (with vertex set $\Z^d$) as follows.  For each $x\in \Z^d$, the directed edge $(x,x+u)$ is in $\mc{G}_x$ if and only if $\omega_x(u)>0$, and the edge set of $\mc{G}(\omega)$ is $\cup_{x\in \Z^d} \mc{G}_x(\omega)$.  For convenience we will also write $\mc{G}=(\mc{G}_x)_{x\in \Z^d}$.    
Note that under $\nu$, $(\mc{G}_x)_{x\in \Z^d}$ are i.i.d.~subsets of $\mc{E}$.  The graph $\mc{G}(\omega)$ is equivalent to the entire graph $\Z^d$, precisely when the environment is {\em elliptic}, i.e.~$\nu(\omega_x(u)>0)=1$ for each $u \in \mc{E}, x\in \Z^d$.  Much of the current literature assumes either the latter condition, or the stronger property of {\em uniform ellipticity}, ie that $\exists \epsilon>0$ such that $\nu(\omega_x(u)>\epsilon)=1$ for each $u \in \mc{E}, x\in \Z^d$.

On the other hand, given a directed graph $\mc{G}=(\mc{G}_x)_{x\in \Z^d}$ (with vertex set $\Z^d$, and such that $\mc{G}_x\ne \varnothing$ for each $x$), we can define a {\em uniform} random environment $\omega=(\omega_x(\mc{G}_x))_{x\in \Z^d}$. Let $|A|$ denote the cardinality of $A$, and set
\[\omega_x(e)=\begin{cases}
|\mc{G}_x|^{-1}, & \text{ if }e \in \mc{G}_x\\
0, & \text{otherwise}
\end{cases}.\]
The corresponding RWRE then moves by choosing uniformly from available steps at its current location.  This gives us a way of constructing rather nice and natural examples of random walks in non-elliptic random environments:  first generate a random directed graph $\mc{G}=(\mc{G}_x)_{x\in \Z^d}$ where $\mc{G}_x$ are i.i.d., then run a random walk on the resulting random graph (choosing uniformly from available steps).  

\begin{DEF}
\label{def:2-val}
We say that the environment is {\em $2$-valued} when $\mu$ charges exactly two points, i.e.~there exist $\gamma_1,\gamma_2\in \mc{P}$ and $p \in (0,1)$ such that $\mu(\{\gamma_1\})=p$, $\mu(\{\gamma_2\})=1-p$.  We say that the graph is {\em $2$-valued} when there exist $E^1,E^2\subset \mc{E}$ and $p\in(0,1)$ such that $\mu(\mc{G}_o=E_1)=p$ and $\mu(\mc{G}_o=E_2)=1-p$.  
\end{DEF}

\cite{HS_RWDRE} proves that the following simple criterion is equivalent to the statement that the random walk visits infinitely many sites almost surely. 

\begin{equation}
\label{unstuckcondition}
\text{There exists an orthogonal set $V$ of unit vectors such that $\mu(\mc{G}_o\cap V\neq\varnothing)=1$.}
\end{equation}

The following is stated in \cite[Lemma 5.1]{HS_RWDRE}:

\begin{LEM} 
\label{lem:trivialtransience}
Assume (\ref{unstuckcondition}) and suppose that $\mu(\downarrow\in \mc{G}_o)=0$ but $\mu(\uparrow\in \mc{G}_o)>0$. 
Then the RWRE is transient in direction $e_2$, $P$-almost surely. 
Let $T$ be the first time the RWRE follows direction $e_2$. If $E[T]<\infty$ then $X_n$ has an asymptotic velocity 
$v=(v^{[1]}, \dots,v^{[d]})$, in the sense that $P(n^{-1}X_n\to v)=1$. Moreover, $v^{[i]}=E[X^{[i]}_T]/E[T].$
\end{LEM}

\proof The random walk visits infinitely many sites, and at each visit to a new site there is positive (non-vanishing) probability of then taking a step in direction $e_2$. Thus the second 
coordinate of the random walk converges monotonically to $\infty$. 

Let $\tau_k$ be the $k$'th time that $X_n$ moves in direction $e_2$, and $\tau_0=0$. Let $Y_k=X_{\tau_k}-X_{\tau_{k-1}}$. Since the environment seen by the random walker is refreshed at every time $\tau_k$, the $Y_k$ are IID, and the $\tau_k$ are sums of IID random variables with distribution that of $T$. Because $E[T]<\infty$, it follows that $E[|Y_k|]<\infty$ as well. By the law of large numbers, $\tau_k/k\to E[T]$ and $X_{\tau_k}/k\to E[Y_1]$ almost surely. Moreover $k^{-1}\max\{|X_n-X_{\tau_{k-1}}|:\tau_{k-1}\le n\le \tau_k\}\to 0$. Thus
$$
\frac{1}{n}X_n\to \frac{1}{E[T]}E[Y_1]=v\quad P\text{-almost surely.}
$$
\Qed\medskip

Table \ref{tab:walks} summarizes what we know about uniform RWDRE in 2-dimensional 2-valued random environments. It reproduces and updates Table 1 of \cite{HS_RWDRE}. 
There is a related table in \cite{HS_DRE1}, giving percolation properties for the directed graphs $\mc{C}$ and $\mc{M}$. The latter includes 2-valued environments such as $(\NSEWalt \,\,\, , \cdot)$ (site percolation), in which one of the possible environments has no arrows. These environments do not appear in the present table, because (as remarked in Section 3 of \cite{HS_RWDRE}), the walk gets stuck on a finite set of vertices (in this case 1 vertex). The RWRE setup we have chosen requires that motion be possible in at least one direction. \medskip

\begin{table}
\begin{center}
\begin{tabular}{l|l|l}
$\gamma_1, \gamma_2$ & Random walk  & Reference  \\
\hline
$\uparrow$ $\rightarrow$ & $v=(1-p,p)$.  &  given here \\
$\uparrow$ $\downarrow$ & Stuck on two vertices.   &  Lemma 2.3 of \cite{HS_RWDRE} \\
$\leftrightarrow$ $\uparrow$ &  $v=\Big(0,\frac{(1-p)^2}{p+(1-p)^2}\Big)$. &  given here \\
$\leftrightarrow$ $\rightarrow$ &  $v=\Big(\frac{1-p}{1+p}, 0\Big)$.  &  given here\\
$\leftrightarrow$ $\updownarrow$ &  $v=(0,0)$. &   $\text{Symmetry}^1$\\
\hline
$\NE$ $\uparrow$ &  $v=\Big(\frac{p}{2}, 1-\frac{p}{2}\Big)$. &  given here\\
$\NE$ $\NW$ &  $v=\Big(\frac{(2p-1)(p^2-p+6)}{6(2-p)(1+p)},\hlf\Big)$.  &   given here\\
$\NE$ $\leftrightarrow$ &  $v=\left(\frac{1}{p^2}+\frac{(1-p)^2}{2p(1-p+p\log p)}\right)^{-1}\cdot (1,1)$.  & given here\\
$\NE$ $\leftarrow$ &  $v=\left(\frac{p(2-p)}{2+3p-2p^2-p^3}\right)\cdot (3,1)+(-1,0)$.  &  given here \\
$\NE$ $\SW$ & $v^{[1]}=v^{[2]}\uparrow$ in $p$. Transient${}^2$ for $p\approx 0,1$.   &   Cor. 2.9 \& Thm. 4.1 of  \cite{HS_RWDRE} 
\\
&{\bf Conjecture:} $v\ne 0$ for $p\neq\hlf$, Recurrent when $p=\hlf$  &\vspace{.1cm} \\
\hline
$\SWE$ $\downarrow$ 
& 
$\frac{1}{v^{[2]}}=\frac{8p(1-p)}{1+\sqrt{5}}-1-2p-\frac{4(1-p)^2(5+\sqrt{5})}{p(1+\sqrt{5})}\overset{\infty}{\underset{n=2}{\sum}}\frac{p^k}{1+2^{-k}(3+\sqrt{5})^k}$, $v^{[1]}=0$.  & given here${}^4$\\
$\SWE$ $\rightarrow$ & $-\frac{1}{v^{[2]}}=4-p-\frac{5+\sqrt{5}}{2}(1-p)^2\Theta(p\gamma)
+\frac{(1-p)[3+\sqrt{5}-(1-p)(5+\sqrt{5})\Theta(p)]^2}{(3+\sqrt{5})[2-(1-p)(5+\sqrt{5})\Theta(p\gamma)]}$ & given here${}^4$\\
& where $\gamma=\frac{3+\sqrt{5}}{2}$ and $\Theta(z)=\sum_{n=0}^\infty \frac{z^n}{\gamma^{2n+1}+1}$. $v^{[1]}=1-3v^{[2]}$. 
 & \\
$\SWE$ $\uparrow$ & $v^{[1]}=0$, $v^{[2]}\downarrow$ in $p$. Transient${}^2$ for $p\approx 0$.  &   Corollary 2.8 of  \cite{HS_RWDRE}  \\
& {\bf Conjecture: } $\exists ! p(\ne 3/4)$ s.t. $v[p]=0$. Recurrent for this $p$. & \\
$\SWE$ $\leftrightarrow$ & $v^{[1]}=0$, $v^{[2]}<0$ for $p>0$. $v^{[2]}$ strictly $\downarrow$ in $p$. & given here${}^5$ \\
$\SWE$ $\updownarrow$ & $v^{[1]}=0$, $v^{[2]}\downarrow$ in $p$. Transient${}^3$ for $p>\frac{3}{4}$,  $v^{[2]}< 0$ for $p>\frac{6}{7}$. & Thm. 4.1 / Thm. 4.10 of \cite{HS_RWDRE}\\
&{\bf Conjecture:} $v^{[2]}< 0$ for $p>0$.   & \\
$\SWE$ $\NE$ & $3v^{[2]}=5v^{[1]}-1$. $v^{[1]}\downarrow$ in $p$. & Thm. 4.1 / Cor. 4.2  of \cite{HS_RWDRE} \\
$\SWE$ $\SW$ & $v^{[1]}=1+3v^{[2]}$ and $v^{[2]}=-1/E[T]$. See below for $E[T]$.   & given here${}^6$\\
$\SWE$ $\NSE$ & $v\cdot(1,-1)=\frac{1}{3}$, $v\cdot (1,1)\downarrow$ in $p$.  & Thm. 4.1 / Cor. 4.2  of \cite{HS_RWDRE} \\
$\SWE$ $\NWE$ &  $v^{[1]}=0$, $v^{[2]}\downarrow$ in $p$. & Thm. 4.1 / Cor. 4.2  of \cite{HS_RWDRE}  \\
& {\bf Conjecture:} $v^{[2]}\ne 0$ for $p\neq\frac12$. Recurrent when $p=\frac12$.&   \vspace{.1cm}\\
\hline 
$\NSEWalt$ \vphantom{${A^A}^A$} \hspace{.1cm} $\uparrow$ & $v^{[1]}=0$, $v^{[2]}\downarrow$ in $p$. Transient${}^3$ for $p<\hlf$, $v^{[2]}> 0$ for $p<\frac13$.  & Thm. 4.1 / Thm. 4.10 of \cite{HS_RWDRE}\\
&  {\bf Conjecture:} $v^{[2]}> 0$ for $p<1$. & \\
$\NSEWalt$ \hspace{.1cm} $\NE$ & $v^{[1]}=v^{[2]}\downarrow$ in $p$. Transient${}^3$ for $p<\hlf$, $v^{[1]}> 0$ for $p<\frac13$.  & Thm. 4.1 / Thm. 4.10 of \cite{HS_RWDRE}\\
&  {\bf Conjecture:}  $v^{[1]}> 0$ for $p<1$.  & \\
$\NSEWalt$ \hspace{.1cm} $\WE$ &  $v=(0,0)$  & Symmetry${}^1$.  \\
$\NSEWalt$ \hspace{.1cm} $\SWE$ & $v^{[1]}=0$, $v^{[2]}\uparrow$ in $p$. 
Transient${}^3$ for $p<\frac{1}{4}$,  $v^{[2]}< 0$ for $p<\frac{1}{7}$. & Thm. 4.1 / Thm. 4.10 of \cite{HS_RWDRE}\\
&{\bf Conjecture:} $v^{[2]}< 0$ for $p<1$.&\\
\hline
\end{tabular}
\end{center}
\caption{Table of results for uniform RWDRE in 2-dimensional 2-valued degenerate random environments, where the first configuration occurs with probability $p\in (0,1)$ and the other with probability $1-p$.}
\label{tab:walks}
\end{table}

\medskip

\noindent {\bf Notes to Table \ref{tab:walks}} \newline
${}^1$ It follows from results of Berger \& Deuschel \cite{BD11} that $\mc{M}$ is recurrent $\forall p\in(0,1)$. \newline
${}^2$ Bounds on the
critical probability are given in \cite{HS_DRE1}. Improved bounds are in preparation.\newline
${}^3$ Improved ranges of values giving transience and speeds are in preparation. \newline
${}^4$ An expansion in terms of $q$-hypergeometric functions is described below.\newline
${}^5$ We do not have a closed form expression for this. But asymptotic expressions are given below.\newline
${}^6$ The expansion is as in the case $(\SWE\rightarrow)$ but messier, so the formula is not included in the table.

\section{Speeds}

\label{sec:speedappendix}

\noindent The non-trivial 2-valued uniform models, in which one must turn to the results of \cite{HS_RWDRE} for existence of a speed, and in which we can at present say very little about the speed (other than mononicity) are as follows:
\begin{itemize}
\item $\NE$ $\SW$\smallskip
\item $\SWE$ $\uparrow$
\item $\SWE$ $\updownarrow$
\item $\SWE$ $\NE$ 
\item $\SWE$ $\NSE$
\item $\SWE$ $\NWE$
\item$\NSEWalt$ \hspace{.1cm} $\uparrow$ 
\item $\NSEWalt$ \hspace{.1cm} $\NE$ \smallskip
\item $\NSEWalt$ \hspace{.1cm} $\SWE$
\end{itemize}

\noindent There are two further models, which are also non-trivial, but for which, once we know that the velocity exists, it must be $v=(0,0)$ by symmetry, namely:
\begin{itemize}
\item $\leftrightarrow\updownarrow$ 
\item $\NSEWalt$ \hspace{.1cm} $\WE$ 
\end{itemize}
 
\noindent The simplest models where one can explicitly calculate the speed are:
 
\begin{itemize}

\item $\uparrow\rightarrow$:  
\newline Because the RWDRE sees a new environment every time, the velocity is simply $(p,1-p)$.

\item $\leftrightarrow\rightarrow$ : 
\newline Let $\tau_k$ be the $k$'th time $n$ that $\mc{G}_{X_n}=\rightarrow$, with $\tau_0$=0. Let $\eta_k=X^{[1]}_{\tau_k}$. At each time $\tau_k$ the process starts exploring a new independent environment, so $T_k=\tau_k-\tau_{k-1}$ are IID (for $k\ge 2$), as are $M_k=\eta_k-\eta_{k-1}$. By the strong law, $\eta_k/\tau_k\to E[M_2]/E[T_2]$ as $k\to\infty$. If $N_n$ is the last $k$ such that $\tau_k\le n$ then 
$$
\frac{\eta_{N_n}}{\tau_{N_{n+1}}}\le \frac{X^{[1]}_n}{n}\le \frac{\eta_{N_{n+1}}}{\tau_{N_n}}
$$
so that $X_n/n\to E[M_2]/E[T_2]$ as well. 

If $M_2=m$ then $\mP_\mc{G}(T_2)=m^2$, since that is the mean time to reach $m$ of a random walk on $[0,m]$ with reflection at 0. Thus $E[T_2]=E[M_2^2]$. But $M_2$ is geometric, $E[M_2]=\sum_{m=1}^\infty mp^{m-1}(1-p)=1/(1-p)$ and $E[M_2^2]=(1+p)/(1-p)^2$. So $v^{[1]}=(1-p)/(1+p)$. 

\item $\NE\uparrow$: 
\newline In this model also, each step of $X_n$ explores a new environment. So we essentially have a regular random walk, whose step distribution is $\rightarrow$ with probability $p/2$ and $\uparrow$ with probability $1-p/2$. So $v=(p/2, 1-p/2)$. 

\end{itemize}

\noindent
In the remaining examples, we use the setup of Lemma \ref{lem:trivialtransience}. There is a direction $e$ for which the first time $T$ that $X_n$ moves in direction $e$ is a renewal time --  what happens starting at time $T$ is independent of what came before. If $e=\pm e_1$, then $v^{[1]}=\pm1/E[T]$. Then if $Y=X^{[2]}_T$, then $v^{[2]}=E[Y]/E[T]$. With corresponding formulae if $e=\pm e_2$. In the following example we calculate $E[Y]$ to get the speed. \medskip

\begin{itemize}

\item $\NE\NW$: 
\newline For $n\ge 0$, let $\tau_n=\inf\{m\ge 0:X_m^{[2]}=n\}$.  Then  for $i\ge 1$, $T_i=\tau_i-\tau_{i-1}$ are i.i.d.~Geometric$(1/2)$ random variables (with mean $2$), and $Y_i=X^{[1]}_{\tau_i -1}-X^{[1]}_{\tau_{i-1}}$ are i.i.d.~random variables, independent of the $\{T_i\}_{i\ge 1}$.  So $E[T_i]=2$ and $v^{[2]}=1/2$. Let $N_n=\sup\{m\ge 0:\tau_m\le n\}$.  Here $e=\uparrow$, and the first time $T$ that we move upwards is geometric with parameter $1/2$. 

Then almost surely,
\[\frac{Y_n^{[1]}}{n}=\frac{\sum_{i=1}^{N_n}Y_i+\sum_{i=\tau_{N_n}+1}^{n}(X^{[1]}_i-X^{[1]}_{i-1})}{n}=\frac{N_n}{n}\frac{\sum_{i=1}^{N_n}Y_i}{N_n}+\frac{\sum_{i=\tau_{N_n}+1}^{n}(X^{[1]}_i-X^{[1]}_{i-1})}{n}\ra \frac{E[Y_1]}{E[T_1]},\]
as $n\ra \infty$, where we have used the fact that $|\sum_{i=\tau_{N_n}+1}^{n}(X^{[1]}_i-X^{[1]}_{i-1})|\le T_{N_{n+1}}$.  

Now let $Y=X^{[1]}_T$, so $v^{[1]}=E[Y]/2$. For $j\ge 1$, we can have $Y=j$ three ways -- reaching no $\NW$ vertex, reaching a $\NW$ vertex at $(j,0)$, or reaching a $\NW$ vertex at $(j+1,0)$. Thus 
\begin{align*}
P(Y=j)&=p^{j+1}\big(\frac12\big)^{j+1}+p^j(1-p)\sum_{n=0}^\infty\big(\frac12\big)^{j+2n+1}+p^{j+1}(1-p)\sum_{n=0}^\infty\big(\frac12\big)^{j+2n+3}\\
&=\frac{p^j(4-p^2)}{3\cdot 2^{j+1}}.
\end{align*}
Likewise, we can have $Y=-j$, $j\ge 1$ three ways, depending on where if anywhere $X_n$ reaches a $\NE$ vertex, giving 
$P(Y=-j)=\big((1-p)^j(4-(1-p)^2)\big)/\big(3\cdot 2^{j+1}\big).$ The case $j=0$ would be similar, but is not needed. Summing over $j$ gives that
\begin{align*}
E[Y]&=\frac{p(4-p^2)}{12}\cdot \frac{1}{(1-p/2)^2}-\frac{(1-p)(4-(1-p)^2)}{12}\cdot\frac{1}{(1-(1-p)/2)^2}\\
&=\frac{p(2+p)}{3(2-p)}-\frac{(1-p)(3-p)}{3(1+p)}
=\frac{(2p-1)(p^2-p+6)}{3(2-p)(1+p)}.
\end{align*}
Comparing this with the speed $\tilde v^{[1]}=p-\hlf$ of a true random walk that goes up with probability $1/2$, right with probability $p/2$, and left with probability $(1-p)/2$, we see that the speeds agree for $p=0, 1/2$, and 1, but the RWRE is slower in between.
See Figure \ref{fig:introcase}.

\begin{figure}
\centering
\includegraphics[scale=.45]{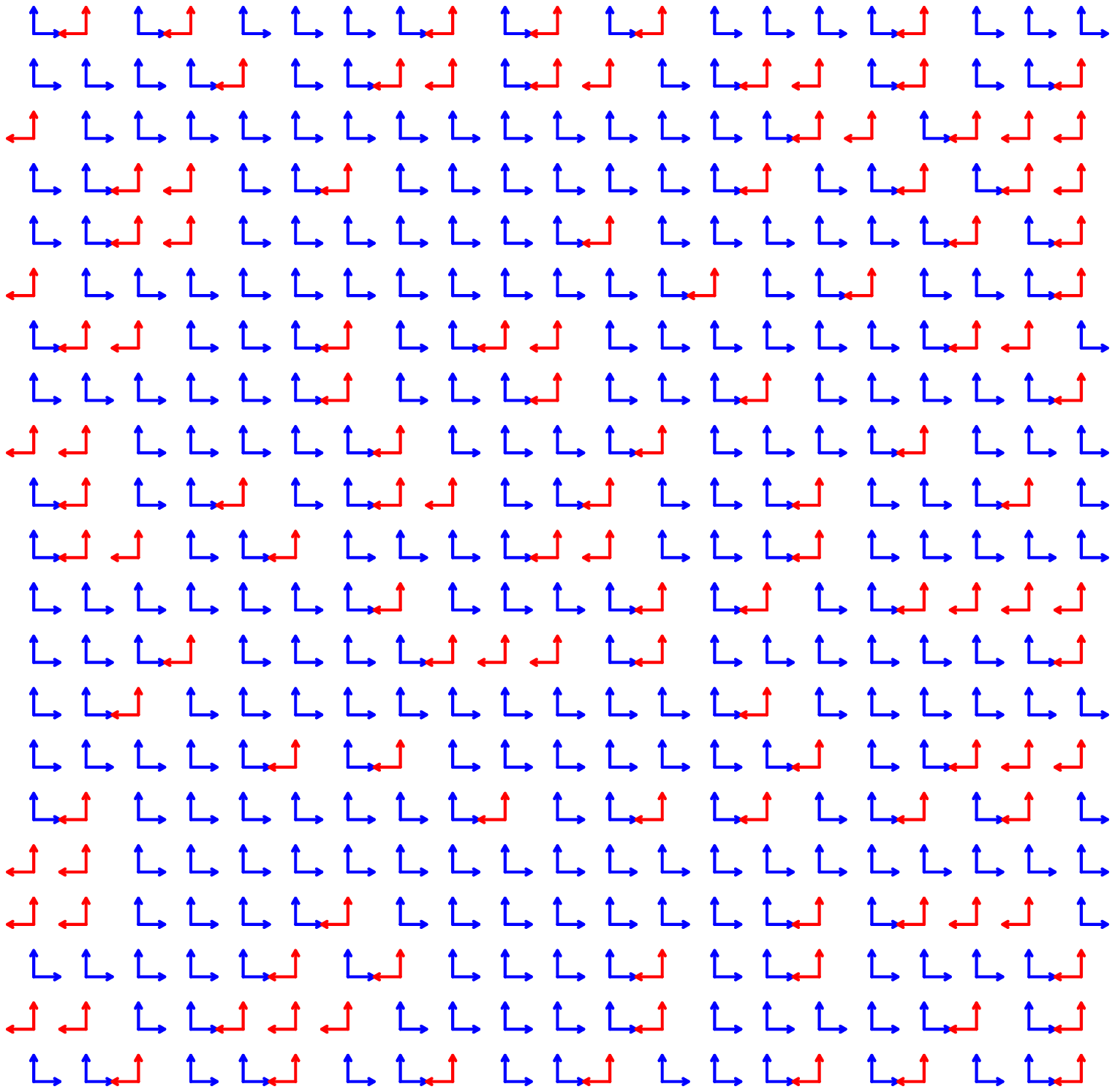}
\includegraphics[height=8cm,width=8cm]{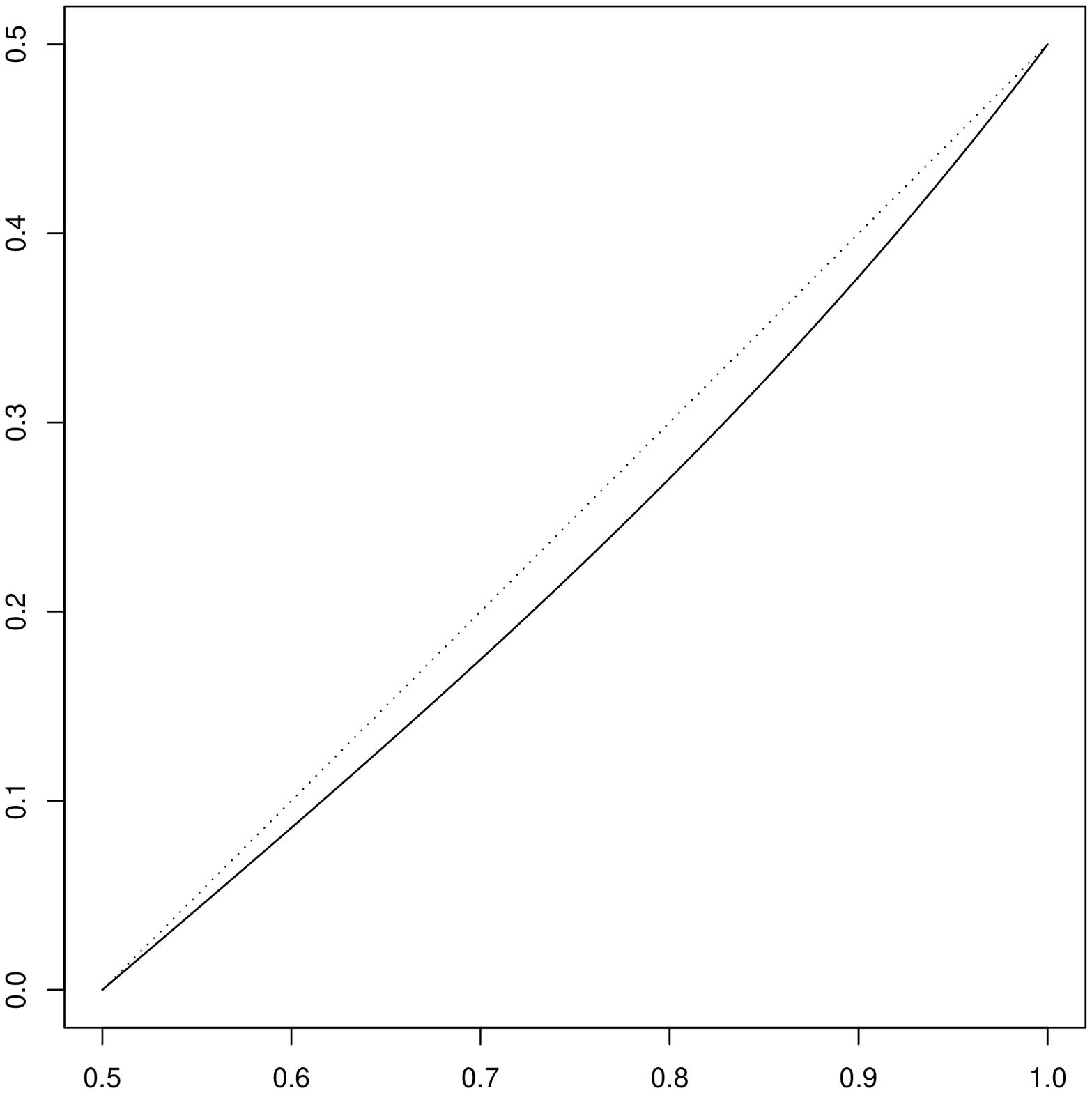}
\caption{A finite region of a degenerate environment in two dimensions such that $\mu(\{\uparrow,\rightarrow\})=p=.75$, $\mu(\{\leftarrow,\uparrow\})=1-p=.25$, and the first coordinate of the velocity as a function of $p\ge \hlf$.}
\label{fig:introcase}
\end{figure}

\end{itemize}

\noindent
In some cases, we can avoid calculating $E[Y]$ directly. Again, assume $e=\pm e_1$. There are two generators, $\mc{L}_1$ and $\mc{L}_2$, depending on the environment. If we apply them to the functions $f_j(x)=x^{[j]}$ we get that
$X_n^{[j]}-\sum_{k<n}(\mc{L}_1f_j1_{\{\mc{G}_{X_k}=\gamma_1\}}+\mc{L}_2f_j1_{\{\mc{G}_{X_k}=\gamma_2\}})$ is a martingale. Therefore
$$
v^{[j]}E[T]=E[X_T^{[j]}]=\mc{L}_1f_j\alpha_1 + \mc{L}_2f_j \alpha_2
$$
where $\alpha_j=E[\#\{k<T:\mc{G}_{X_k}=A_j\}]$. When $j=1$ the LHS is $\pm 1$, which usually lets us solve for $\alpha_1$ in terms of $E[T]$. We know that $\alpha_2=E[T]-\alpha_1$, so putting $j=2$ will then give us $v^{[2]}$. Thus all that remains is to calculate $E[T]$. (We could have done the previous example this was as well.)

\begin{itemize}

\item $\leftrightarrow\uparrow$ : 
\newline Here $e=\uparrow$. We have $v^{[1]}=0$ by symmetry, and $v^{[2]}=1/E[T]$. If the origin is $\uparrow$, then $T=1$. Otherwise, suppose there are $\uparrow$ at $(-i,0)$ and at $(j,0)$, with only $\WE$ in between. Then $\mP_{\mc{G}}(T-1)=ij$, since this is the mean exit time for a simple random walk on $[-i,j]$. Given that the origin is $\WE$ (which happens with probability $p$), $i$ and $j$ are independent geometric random variables, with means $1/(1-p)$. Thus $E[T]=1+p/(1-p)^2$.\medskip

\item $\NE\leftrightarrow$: 
\newline Here $e=\uparrow$, and the martingale equations are that $1=v^{[2]}E[T]=\alpha_2/2$ and $v^{[1]}E[T]=\alpha_2/2$. In other words, $v^{[1]}=v^{[2]}=1/E[T]$. So we must now find $E[T]$. 

First consider a random walk $Z_j$ on 
$[0,n)$ with the following boundary conditions: at $n$ there is absorption, and at $0$ we reflect with probability $1/2$ and die otherwise. Let $S$ be the time of death or absorption, and let $f(k)=E[S\mid Z_0=k]$. Then 
$$
f(k)=1+\frac{f(k-1)+f(k+1)}{2}
$$
for $1\le k\le n-1$, $f(n)=0$, and $f(0)=1+f(1)/2$. The solution to the recurrence is $f(k)=A+Bk-k^2$, and substituting the boundary conditions gives $f(k)=(n-k)(k+1)$. Likewise let $g(k)=P(Z_S=n\mid Z_0=k)$. Then $g(k)=[g(k-1)+g(k+1)]/2$ for $1\le k\le n-1$, with boundary conditions $g(0)=g(1)/2$ and $g(n)=1$. This has solution $g(k)=(k+1)/(n+1)$. 

Now think of how $X_j$ evolves. Let the first $\NE$ to the left of $o$ be at $x_0=(i_0,0)$, where $i_0\le 0$. Let successive $\NE$ to the right of $o$ be at $x_1=(i_1,0)$, $x_2=(i_2,0)$, etc., where $0<i_1<i_2<\dots$. On the horizontal interval 
$[x_0,x_1)$, $X_j$ performs a simple random walk till it hits $x_0$ or $x_1$. When it hits $x_0$ it either moves upward (making this time $T$), or it reflects back into the interval. If it reaches $x_1$ it leaves this interval forever, and starts the same process over again on the interval 
$[x_1,x_2)$. Let the interval being visited at time $T$ be 
$[x_N,x_{N+1})$, where $N\ge 0$, and let $S_j$ be the total time spent in 
$[x_j, x_{j+1})$. Therefore 
$T=\sum_{j=0}^NS_j$, and 
$$
E[T]=\sum_{j=0}^\infty E[S_j 1_{\{N\ge j\}}].
$$
Let $A_j$ be the event that $X_n$ exits 
$[x_j, x_{j+1})$ on the right, ie. at $x_{j+1}$. Then for $j\ge 1$, $\{N\ge j\}=\cap_{k=0}^{j-1}A_k$. Moreover, there is a renewal every time $X_n$ enters a new interval, because a new environment starts getting explored. The interval 
$[x_0,x_1)$ is different from the rest, because we start at $o$. But for all other 
$[x_j,x_{j+1})$ the process starts walking at $x_j$. Therefore the cases $j\ge 1$ are actually independent replications of the same procedure. In other words,
$$
E[T]=E[S_0]+\sum_{j=1}^\infty E[S_j \prod_{k=0}^{j-1}1_{A_k}]=E[S_0]+\sum_{j=1}^\infty E[S_1\mid A_0] P(A_0)P(A_1\mid A_0)^{j-1}.
$$
We use the expressions for $f$ and $g$ to work out these factors. By the expression for $f$, we have $\mP_{\mc{G}}(S_0)=i_1(1-i_0)$. Take $i=-i_0$ and $k=i_1$. Then
$$
E[S_0]=\sum_{i=0}^\infty\sum_{k=1}^\infty p^2(1-p)^{i+k-1}k(1+i)
=\Big(\sum_{k=1}^\infty p(1-p)^{k-1}k\Big)^2=\frac{1}{p^2}.
$$
Likewise $\mP_{\mc{G}}(S_1\mid A_0)=i_2-i_1$. Write $k$ for this quantity, so 
$$
E[S_1\mid A_0]=\sum_{k=1}^\infty (1-p)^{k-1}pk=\frac{p}{p^2}=\frac{1}{p}.
$$
Similarly, $\mP_{\mc{G}}(A_0)=(1-i_0)/(1+i_1-i_0)$, so (letting $n=i+k$)
\begin{align*}
P(A_0)&=\sum_{i=0}^\infty\sum_{k=1}^\infty p^2(1-p)^{i+k-1}\frac{1+i}{1+i+k}
=\sum_{n=1}^\infty p^2(1-p)^{n-1}\sum_{j=1}^k\frac{n+1-j}{1+n}\\
&=\sum_{n=1}^\infty \frac{p^2(1-p)^{n-1}}{n+1}\cdot\frac{n(n+1)}{2}
=\frac{p^2}{2}\cdot\frac{1}{p^2}=\frac12.
\end{align*}
And $\mP_{\mc{G}}(A_1\mid A_0)=1/(1+i_2-i_1)$, so 
\begin{align*}
P(A_1\mid A_0)&=\sum_{k=1}^\infty p(1-p)^{k-1}\frac{1}{1+k}
=\frac{p}{(1-p)^2}\Big(\sum_{k=0}^\infty \frac{(1-p)^{n+1}}{n+1}-(1-p)\Big)\\
&=\frac{p}{(1-p)^2}(-\log p-(1-p))
=1-\frac{1-p+p\log p}{(1-p)^2}.
\end{align*}
Putting this together, 
$$
E[T]=\frac{1}{p^2}+\frac{(1-p)^2}{2p(1-p+p\log p)}.
$$

\item $\NE\leftarrow$ :
\newline Here $e=\uparrow$, and the martingale equations are that $1=v^{[2]}E[T]=\alpha_1/2$ and $v^{[1]}E[T]=\alpha_1/2-\alpha_2=-E[T]+3\alpha_1/2$. In other words, $v^{[2]}=1/E[T]$ and $v^{[1]}=-1+3/E[T]$. So we must now find $E[T]$. 

Suppose that there is a $\NE$ at $(-i,0)$ for $i\ge 1$, and $\leftarrow$'s at $o$ and all points in between (a scenario with probability $p(1-p)^i$. Then $X_n$ takes $i$ steps to the left, and then oscillates between $(-i,0)$ and $(-i+1,0)$ a random number of times, before $T$ occurs. 

The other possibility is that there is a $\leftarrow$ at $(j,0)$ for $j\ge 1$, and $\NE$'s at $o$ and all points in between. This scenario has probability $(1-p)p^j$. Now $X_n$ steps right, and $T$ may occur before it reaches $(j),j$, or it may reach $(j,0)$ and then oscillate until time $T$. The various scenarios lead to the following expression:
\begin{align*}
E[T]&=\sum_{i=1}^\infty p(1-p)^i\sum_{k=0}^\infty (1/2)^{k+1}[i+2k+1]\\
&\qquad\qquad+\sum_{j=1}^\infty (1-p)p^j\Big(\sum_{k=0}^{j-2}(1/2)^{k+1}[k+1] + \sum_{k=0}^\infty (1/2)^{j+k}[j+2k]\Big)\\
&=\sum_{i=1}^\infty p(1-p)^i(i+3)+\sum_{j=1}^\infty (1-p)p^j\Big(\frac12\frac{d}{dt}\Big|_{t=1/2}\frac{1-t^j}{1-t}+\frac{j+2}{2^{j-1}}\Big)\\
&=3(1-p)+\frac{1-p}{p}+2\sum_{j=1}^\infty (1-p)p^j\Big(1+\frac{1}{2^j}\Big)\\
&=2-3p+\frac{1}{p}+2p+\frac{p(1-p)}{1-p/2}
=\frac{2+3p-2p^2-p^3}{p(2-p)}.
\end{align*}

\item $\SWE\downarrow $:
\newline The velocity is $(0,-1/E[T])$, where $T$ is the time of the first step in the $\downarrow$ direction. 

To find $E[T]$, first consider random walk on $[0,n]$, with the probability of death in 1 step starting from $1\le k\le n-1$ being $1/3$, and the probability of death in 1 step being 1 starting from 0 or $n$. Let $f(k)$ be the mean time of death, starting from $k$. Then $f(0)=f(n)=1$, and otherwise
$$
f(k)=1+\frac{f(k-1)+f(k+1)}{3}.
$$
The solution is $f(k)=3+\bar c\bar\gamma^k + c\gamma^k$ where $\bar \gamma<\gamma$ are solutions of $z +z^{-1}=3$. In other words, $\gamma=(3+\sqrt{5})/2$ and $\bar \gamma=(3-\sqrt{5})/2$ . From the boundary conditions, we get
$$
f(k)=3+\frac{2(1-\gamma^n)\bar\gamma^k}{\gamma^n-\bar\gamma^n}+\frac{2(1-\bar\gamma^n)\gamma^k}{\gamma^n-\bar\gamma^n}.
$$
But $\bar\gamma\gamma=1$, so this simplifies to 
$$
f(k)=3-2\frac{\gamma^k+\gamma^{n-k}}{\gamma^n+1}.
$$
If there is a $\downarrow$ at $o$ then $T=1$. Otherwise, suppose there are $\downarrow$'s at $(-i,0)$ and $(j,0)$, with $\SWE$'s in between, where $i,j\ge 1$. Then $\mP_{\mc{G}}(T)=3-2(\gamma^i+\gamma^j)/(\gamma^{i+j}+1)$. Therefore 
\begin{align*}
E[T]
&=(1-p)\cdot 1 +p\sum_{i,j=1}^\infty p^{i+j-2}(1-p)^2\Big(3-2\frac{\gamma^i+\gamma^j}{\gamma^{i+j}+1}\Big)\\
&=1-p+3p(1-p)^2\Big(\sum_{i=1}^\infty p^{i-1}\Big)^2 -\frac{2(1-p)^2}{p}\sum_{k=2}^\infty \frac{p^k}{\gamma^k+1}\sum_{j=1}^{k-1}(\gamma^j+\gamma^{k-j})\\
&=1+2p-\frac{4(1-p)^2}{p}\sum_{k=2}^\infty \frac{p^k}{\gamma^k+1}\Big(\frac{\gamma^k-1}{\gamma-1}-1\Big)\\
&=1+2p-\frac{4(1-p)^2}{p}\sum_{k=2}^\infty \Big(\frac{1}{\gamma-1}p^k-\frac{\gamma+1}{\gamma-1}\frac{p^k}{\gamma^k+1}\Big)\\
&=1+2p-\frac{4p(1-p)}{\gamma-1}+\frac{4(1-p)^2(\gamma+1)}{p(\gamma-1)}\sum_{k=2}^\infty \frac{p^k}{\gamma^k+1}.
\end{align*}
Note that the expression $Q(q;p)=1+2\sum_{k=1}^\infty \frac{p^k}{q^k+1}$ is known as a unilateral $q$-hypergeometric series, and in the theory of special functions would be written
$$
Q(q;p)={}_2\phi_1\left[
\begin{matrix}
q & -1\\
-q & 
\end{matrix}\,; q, p\right].
$$
See \cite{Wiki}. We are indebted to Martin Muldoon for pointing this out.
We could therefore also write 
\begin{align*}
E[T]&=\frac{2(1-p)^2(\gamma+1)}{p(\gamma-1)}Q(\gamma;p)+1+2p\\
&\qquad\qquad\qquad-\frac{4p(1-p)}{\gamma-1}-\frac{4(1-p)^2}{\gamma-1}-\frac{2(1-p)^2(\gamma+1)}{p(\gamma-1)}\\
&=\frac{2(1-p)^2(\gamma+1)}{p(\gamma-1)}Q(\gamma;p)+5-\frac{2}{p}-\frac{4(1-p)}{p(\gamma-1)}
\end{align*}

\item 
$\SWE\rightarrow $: 
\newline Here $e=\downarrow$, and the martingale equations are that $-1=v^{[2]}E[T]=-\alpha_1/3$ and $v^{[1]}E[T]=\alpha_2=E[T]-\alpha_1=E[T]-3$. In other words, $v^{[2]}=-1/E[T]$ and $v^{[1]}=1-3/E[T]=1+3v^{[2]}$. So we must now find $E[T]$. 

First consider a random walk $Z_j$ on $[0,n]$ with the following boundary conditions: it reflects at $0$, and it is absorbed at $n$.  At points in between there is killing with probability $1/3$, and otherwise $Z_j$ performs a simple symmetric random walk. Let $S$ be the time of death or absorption, and let $f(k)=f_n(k)=E[S\mid Z_0=k]$. Then 
$$
f(k)=1+\frac{f(k-1)+f(k+1)}{3}
$$
for $1\le k\le n-1$, $f(n)=0$, and $f(0)=1+f(1)$. The solution to the recurrence is $f(k)=3+\bar c\bar\gamma^k+c\gamma^k$ where as above, $\bar\gamma<\gamma$ are  
$(3\pm\sqrt{5})/2$. From the boundary conditions, we get
$$
f(k)=3+\bar\gamma^k\frac{3\gamma-3-\gamma^n}{\gamma^n(\bar\gamma-1)-\bar\gamma^n(\gamma-1)}+\gamma^k\frac{\bar\gamma^n-3\bar\gamma+3}{\gamma^n(\bar\gamma-1)-\bar\gamma^n(\gamma-1)}.
$$
But $\bar\gamma\gamma=1$, so this simplifies to 
\begin{align*}
f(k)&=3-\frac{\gamma^{n-k}(3\gamma-3-\gamma^n)+\gamma^k(1-3\gamma^{n-1}+3\gamma^n)}{(\gamma-1)(\gamma^{2n-1}+1)}\\
&=3-\frac{1}{\gamma^{2n-1}+1}\Big[3(\gamma^{n-k}+\gamma^{n+k-1})+\frac{\gamma^k-\gamma^{2n-k}}{\gamma-1}\Big].
\end{align*}
Likewise let $g(k)=g_n(k)=P(Z_S=n\mid Z_0=k)$. Then $g(k)=[g(k-1)+g(k+1)]/3$ for $1\le k\le n-1$, so $g(k)=\bar{c}\bar{\gamma}^k+c\gamma^k$ with boundary conditions $g(0)=g(1)$ and $g(n)=1$. As above, this has solution 
$$
g(k)=\frac{\gamma^{n-k}+\gamma^{n+k-1}}{\gamma^{2n-1}+1}.
$$

Now consider how $X_j$ evolves. Let the first $\rightarrow$ to the left of $o$ be at $x_0=(i_0,0)$, where $i_0\le 0$. Let successive $\rightarrow$ to the right of $o$ be at $x_1=(i_1,0)$, $x_2=(i_2,0)$, etc., where $0<i_1<i_2<\dots$. At interior points of the horizontal interval 
$[x_0,x_1)$, $X_j$ leaves the interval with probability $1/3$ at every step, by moving downwards (ie at time $T$). Otherwise it evolves as a simple random walk till it hits $x_0$ or $x_1$. If it hits $x_0$ it reflects back with probability 1. If it reaches $x_1$ it leaves this interval forever, and starts the same process over again on the interval 
$[x_1,x_2)$. Let the interval being visited at time $T$ be 
$[x_N,x_{N+1})$, where $N\ge 0$, and let $S_j$ be the total time spent in 
$[x_j, x_{j+1})$. Then 
$T=\sum_{j=0}^NS_j$, and 
$$
E[T]=\sum_{j=0}^\infty E[S_j 1_{\{N\ge j\}}].
$$
If $A_j$ is the event that $X_n$ exits 
$[x_j, x_{j+1})$ at $x_{j+1}$ (ie before $T$), then for $j\ge 1$ we have $\{N\ge j\}=\cap_{k=0}^{j-1}A_k$. Moreover, there is a renewal every time $X_n$ enters a new interval, because we start exploring a new environment. The interval 
$[x_0,x_1)$ is different from the rest, because we start at $o$. But for all other 
$[x_j,x_{j+1})$ the process starts walking at $x_j$. In other words, the cases $j\ge 1$ are independent replications of the same procedure. Therefore
\begin{align*}
E[T]&=E[S_0]+\sum_{j=1}^\infty E[S_j \prod_{k=0}^{j-1}1_{A_k}]=E[S_0]+\sum_{j=1}^\infty E[S_1\mid A_0] P(A_0)P(A_1\mid A_0)^{j-1}\\
&=E[S_0]+E[S_1\mid A_0] P(A_0)/(1-P(A_1\mid A_0)).
\end{align*}
We can work out all these factors using the expressions for $f$ and $g$. First define
$$
\Theta(z)=\Theta_\gamma(z)=\sum_{n=0}^\infty\frac{z^n}{\gamma^{2n+1}+1}=\frac{1}{2\sqrt{z}}\Big[Q(\gamma;\sqrt{z})-Q(\gamma^2;z)\Big]
$$
where $Q$ is the $q$-hypergeometric series defined earlier. 
Then using that $\gamma^2-3\gamma+1=0$, 
\begin{align*}
E[S_0]&=\sum_{i_0=0}^\infty \sum_{i_1=1}^\infty p^{i_0+i_1-1}(1-p)^2 f_{i_0+i_1}(i_0)\\
&=3-(1-p)^2\sum_{n=1}^\infty p^{n-1}\sum_{k=0}^{n-1}\Big[\frac{3(\gamma^{n-k}+\gamma^{n+k-1})}{\gamma^{2n-1}+1}
+\frac{\gamma^k-\gamma^{2n-k}}{(\gamma-1)(\gamma^{2n-1}+1)}\Big]\\
&=3-(1-p)^2\sum_{n=1}^\infty\frac{p^{n-1}(\gamma^n-1)}{(\gamma-1)^2(\gamma^{2n-1}+1)}\Big[3(\gamma-1)(\gamma+\gamma^{n-1})+1-\gamma^{n+1}\Big]\\
&=3-\frac{(1-p)^2}{(\gamma-1)^2}\sum_{n=1}^\infty\frac{p^{n-1}}{\gamma^{2n-1}+1}\Big[(\gamma^{2n-1}+1)(3\gamma-3-\gamma^2)+\gamma^{n-1}(3\gamma^3-2\gamma^2-2\gamma+3)-2(\gamma^2-1)\Big]\\
&=3-\frac{(1-p)(3\gamma-3-\gamma^2)}{(\gamma-1)^2}+\left(\frac{1-p}{\gamma-1}\right)^2\Big[2(\gamma^2-1)\Theta(p)-(3\gamma^3-2\gamma^2-2\gamma+3)\Theta(p\gamma)\Big],\\
&=3+\frac{2(1-p)}{\gamma}+\frac{2(1-p)^2}{\gamma}\Big[(3\gamma-2)\Theta(p)-(8\gamma-2)\Theta(p\gamma)\Big],\\
E[S_1\mid A_0]&=\sum_{n=1}^\infty p^{n-1}(1-p)f_n(0)
=3-\sum_{n=1}^\infty \frac{p^{n-1}(1-p)}{\gamma^{2n-1}+1}\Big[3\gamma^{n-1}(1+\gamma)-\frac{\gamma^{2n}-1}{\gamma-1}\Big]\\
&=3-3(1-p)(1+\gamma)\Theta(p\gamma)+\frac{\gamma}{\gamma-1}-\frac{(1-p)(1+\gamma)}{\gamma-1}\Theta(p),
\\
P(A_0)&=\sum_{i_0=0}^\infty \sum_{i_1=1}^\infty p^{i_0+i_1-1}(1-p)^2 g_{i_0+i_1}(i_0)
=(1-p)^2\sum_{n=1}^\infty p^{n-1}\sum_{k=0}^{n-1}\frac{\gamma^{n-k}+\gamma^{n+k-1}}{\gamma^{2n-1}+1}\\
&=\frac{(1-p)^2}{\gamma-1}\sum_{n=1}^\infty \frac{p^{n-1}(\gamma^{2n-1}+\gamma^{n+1}-\gamma^{n-1}-\gamma)}{\gamma^{2n-1}+1}
=\frac{(1-p)}{\gamma-1}+(1-p)^2(\gamma+1)\Big[\Theta(p\gamma)-\frac{\Theta(p)}{\gamma-1}\Big],\\
P(A_1\mid A_0)&=\sum_{n=1}^\infty p^{n-1}(1-p)g_n(0)=(1-p)(\gamma+1)\Theta(p\gamma).
\end{align*}
Therefore 
\begin{multline*}
E[T]=3+\frac{2(1-p)}{\gamma}+\frac{2(1-p)^2}{\gamma}\Big[(3\gamma-2)\Theta(p)-(8\gamma-2)\Theta(p\gamma)\Big]\\
+\frac{\Big[3-3(1-p)(1+\gamma)\Theta(p\gamma)+\frac{\gamma}{\gamma-1}-\frac{(1-p)(\gamma+1)}{\gamma-1}\Theta(p)\Big]\Big[\frac{(1-p)}{\gamma-1}+(1-p)^2(\gamma+1)\Big\{\Theta(p\gamma)-\frac{\Theta(p)}{\gamma-1}\Big\}\Big]}{1-(1-p)(\gamma+1)\Theta(p\gamma)}.
\end{multline*}
Simplifying this, we have
$$
E[T]=4-p-(1-p)^2\frac{4\gamma-1}{\gamma}\Theta(p\gamma)
+\frac{1-p}{\gamma}\frac{\Big[\gamma-(1-p)(\gamma+1)\Theta(p)\Big]^2}{1-(1-p)(\gamma+1)\Theta(p\gamma)}.
$$
Substituting for $\gamma$ gives
$$
E[T]=4-p-\frac{5+\sqrt{5}}{2}(1-p)^2\Theta(p\gamma)
+\frac{(1-p)\Big[3+\sqrt{5}-(1-p)(5+\sqrt{5})\Theta(p)\Big]^2}{(3+\sqrt{5})\Big[2-(1-p)(5+\sqrt{5})\Theta(p\gamma)\Big]}.
$$
\smallskip

\item $\SWE\SW $:
\newline Here $e=\downarrow$, and the martingale equations are that $-1=v^{[2]}E[T]=-\alpha_1/3-\alpha_2/2=-E[T]/3-\alpha_2/6$ and $v^{[1]}E[T]=-\alpha_2/2=E[T]-3$. In other words, $v^{[2]}=-1/E[T]$ and $v^{[1]}=1-3/E[T]=1+3v^{[2]}$. So we must now find $E[T]$. This is just like the previous example, but messier.\medskip

It is slightly more convenient to work out $E[T]$ using $(\SWE\SE)$, since the structure of the previous example can be maintained with minor changes. For $f$ and $g$, the difference is that now there is killing at 0 at rate $1/2$. This changes the boundary condition for $f(0)$ to be 
$f(0)=1+\frac12 f(1)$, from which we get
$$
f(k)=3+\bar{\gamma}^k\frac{\gamma^n+3\gamma-6}{\gamma^n(\bar{\gamma}-2)-\bar{\gamma}^n(\gamma-2)}+\gamma^k\frac{6-\bar{\gamma}^n-3\bar{\gamma}}{\gamma^n(\bar{\gamma}-2)-\bar{\gamma}^n(\gamma-2)}
$$
and therefore 
$$
f(k)=3+\frac{1}{\gamma^{2n-1}(1-2\gamma)-\gamma+2}\Big[ \gamma^{n-k}(\gamma^n+3\gamma-6)+\gamma^k(6\gamma^n-3\gamma^{n-1}-1)\Big].
$$
The changed boundary condition for $g$ is that $g(0)=\frac12 g(1)$, from which we get
$$
g(k)=\bar{\gamma}^k\frac{2-\gamma}{\gamma^n(\bar{\gamma}-2)-\bar{\gamma}^n(\gamma-2)}+\gamma^k\frac{\bar{\gamma}-2}{\gamma^n(\bar{\gamma}-2)-\bar{\gamma}^n(\gamma-2)}
$$
and therefore
$$
g(k)=\frac{\gamma^{n-k}(2-\gamma)+\gamma^{n+k-1}(1-2\gamma)}{\gamma^{2n-1}(1-2\gamma)-\gamma+2}.
$$
For the walk $X_j$ there is now reflection at $x_k$ with probability $\frac12$ and $\downarrow$ with probability $\frac12$. Set
$$
\tilde\Theta(z)=\tilde \Theta_\gamma(z)=\sum_{n=0}^\infty\frac{z^n}{\gamma^{2n+1}(1-2\gamma)-\gamma+2}.
$$
\begin{align*}
E[S_0]&=\sum_{i_0=0}^\infty \sum_{i_1=1}^\infty p^{i_0+i_1-1}(1-p)^2 f_{i_0+i_1}(i_0)\\
&=3+(1-p)^2\sum_{n=1}^\infty \frac{p^{n-1}}{\gamma^{2n-1}(1-2\gamma)-\gamma+2}\sum_{k=0}^{n-1}\Big[ \gamma^{n-k}(\gamma^n+3\gamma-6)+\gamma^k(6\gamma^n-3\gamma^{n-1}-1)\Big]\\
&=3+\frac{(1-p)^2}{\gamma-1}\sum_{n=1}^\infty \frac{p^{n-1}(\gamma^n-1)}{\gamma^{2n-1}(1-2\gamma)-\gamma+2}\Big[ \gamma(\gamma^n+3\gamma-6)+6\gamma^n-3\gamma^{n-1}-1\Big]\\
&=3+\frac{(1-p)^2}{\gamma-1}\sum_{n=1}^\infty \frac{p^{n-1}}{\gamma^{2n-1}(1-2\gamma)-\gamma+2}\Big[ \gamma^{2(n-1)}(\gamma^3+6\gamma^2-3\gamma)+\\
&\hspace{3in}+\gamma^{n-1}(3\gamma^3-7\gamma^2-7\gamma+3)-(3\gamma^2-6\gamma+1)\Big]\\
&=3+\frac{(1-p)^2}{\gamma-1}\Big[ (\gamma^3+6\gamma^2-3\gamma)\tilde\Theta(p\gamma^2)+(3\gamma^3-7\gamma^2-7\gamma+3)\tilde\Theta(p\gamma)-(3\gamma^2-6\gamma+1)\tilde\Theta(p)\Big]\\
E[S_1\mid A_0]&=\sum_{n=1}^\infty p^{n-1}(1-p)f_n(0)\\
&=3+\sum_{n=1}^\infty \frac{p^{n-1}(1-p)}{\gamma^{2n-1}(1-2\gamma)-\gamma+2}\Big[\gamma^{n}(\gamma^n+3\gamma-6)+(6\gamma^n-3\gamma^{n-1}-1)\Big]\\
&=3+\sum_{n=1}^\infty \frac{p^{n-1}(1-p)}{\gamma^{2n-1}(1-2\gamma)-\gamma+2}\Big[\gamma^2\gamma^{2(n-1)}+3\gamma^{n-1}(\gamma^2-1)-1\Big]\\
&=3+(1-p)\Big[\gamma^2\tilde\Theta(p\gamma^2)+3(\gamma^2-1)\tilde\Theta(p\gamma)-\tilde\Theta(p)\Big]\\
P(A_0)&=\sum_{i_0=0}^\infty \sum_{i_1=1}^\infty p^{i_0+i_1-1}(1-p)^2 g_{i_0+i_1}(i_0)
=(1-p)^2\sum_{n=1}^\infty p^{n-1}\sum_{k=0}^{n-1}\frac{\gamma^{n-k}(2-\gamma)+\gamma^{n+k-1}(1-2\gamma)}{\gamma^{2n-1}(1-2\gamma)-\gamma+2}\\
&=\frac{(1-p)^2}{\gamma-1}\sum_{n=1}^\infty \frac{p^{n-1}}{\gamma^{2n-1}(1-2\gamma)-\gamma+2}\Big[\gamma(\gamma^n-1)(2-\gamma)+\gamma^{n-1}(\gamma^n-1)(1-2\gamma)\Big]\\
&=\frac{(1-p)^2}{\gamma-1}\sum_{n=1}^\infty \frac{p^{n-1}}{\gamma^{2n-1}(1-2\gamma)-\gamma+2}\Big[\gamma^{2(n-1)}\gamma(1-2\gamma)-\gamma^{n-1}(\gamma^3-2\gamma^2-2\gamma+1)+\gamma(\gamma-2)\Big]\\
&=\frac{(1-p)^2}{\gamma-1}\Big[(\gamma(1-2\gamma)\tilde\Theta(p\gamma^2)-(\gamma^3-2\gamma^2-2\gamma+1)\tilde\Theta(p\gamma)+\gamma(\gamma-2)\tilde\Theta(p)\Big]\\
P(A_1\mid A_0)&=\sum_{n=1}^\infty p^{n-1}(1-p)g_n(0)=(1-p)\sum_{n=1}^\infty \frac{p^{n-1}}{\gamma^{2n-1}(1-2\gamma)-\gamma+2}\Big[\gamma^{n-1}(1-\gamma^2)\Big]\\
&=(1-p)(1-\gamma^2)\tilde\Theta(p\gamma).
\end{align*}
Therefore
\begin{multline*}
E[T]=3+\frac{(1-p)^2}{\gamma-1}\Big[ (\gamma^3+6\gamma^2-3\gamma)\tilde\Theta(p\gamma^2)+(3\gamma^3-7\gamma^2-7\gamma+3)\tilde\Theta(p\gamma)-(3\gamma^2-6\gamma+1)\tilde\Theta(p)\Big]+\\
+\Big[3+(1-p)\big[\gamma^2\tilde\Theta(p\gamma^2)+3(\gamma^2-1)\tilde\Theta(p\gamma)-\tilde\Theta(p)\big]\Big]\frac{(1-p)^2}{\gamma-1}\cdot\\
\cdot\Big[(\gamma(1-2\gamma)\tilde\Theta(p\gamma^2)-(\gamma^3-2\gamma^2-2\gamma+1)\tilde\Theta(p\gamma)+\gamma(\gamma-2)\tilde\Theta(p)\Big]/\Big[1-(1-p)(1-\gamma^2)\tilde\Theta(p\gamma)\Big]
\end{multline*}

\end{itemize}

\noindent
Finally, we give an asymptotic argument, in the one case in which there is an elementary renewal structure for which we don't know how to find the speed analytically. 

\begin{itemize}

\item $\SWE\WE $: 
\newline Here $e=\downarrow$, and $v=(0,-1/E[T])$. Though we don't know how to find $v$ analytically, here is an approach that should give an asymptotic expansion in powers of $q=1-p$. 

Embed $\Z\subset\Z^2$ as $\Z\times\{0\}$, and let $Y_n=X^{[1]}_n$ for $n<T$. We can fill in new independent increments after time $T-1$ to make $Y_n$ into a simple symmetric random walk started from $0$, and then recover an independent copy of $T$ by killing $Y$ at a rate that depends on the environment. Write $\tilde P$ for this extension of $P$. Let $V_i$ be 1 (resp. $2/3$) if $\mc{G}_{(i,0)}=\WE$ (resp. $\SWE$). Then by Feynman-Kac, 
$$
E[T]=\sum_{k=1}^\infty P(T\ge k)=\sum_{k=1}^\infty \tilde P(\prod_{i=0}^{k-2} V_{Y_i}).
$$
This equals 
$\sum_{k=1}^\infty \tilde P(\prod_{j\in\Z} V_j^{N_j(k-1)})$, where $N_j(\ell)$ counts the number of visits of $Y_n$ to $j$, for $0\le n<\ell$. Integrating out the environment, and setting $\gamma=2/3$, we get
\begin{align*}
\sum_{k=1}^\infty \tilde P(\prod_{j\in\Z} [1-p+p\gamma^{N_j(k-1)}])
&= \sum_{k=0}^\infty \tilde P(\prod_{j\in\Z} \gamma^{N_j(k)}[q\gamma^{-N_j(k)}+1-q])\\
&= \sum_{k=0}^\infty \gamma^{k}\tilde P(\prod_{j\in\Z} [1+q(\gamma^{-N_j(k)}-1)])
\end{align*}
(since $\sum_j N_j(k)=k$). This expression can in principle be expanded as a series in $q$. What complicates this is that the $q^i$ term involves knowing the joint distributions of the $N_j(k)$ for $i$ choices of $j$. Still, it is easy to work out the constant and linear terms (and with some work one could in principal take this further). In particular, the following should be true:
\begin{align*}
E[T]&=\sum_{k=0}^\infty \gamma^{k} + q\sum_{k=0}^\infty \sum_{j\in\Z} \gamma^{k} \tilde P(\gamma^{-N_j(k)}-1)+O(q^2)\\
&=3+q \sum_{j\in\Z} \Big[\sum_{k=0}^\infty \gamma^{k} \tilde P(\gamma^{-N_j(k)})-3\Big]+O(q^2).
\end{align*}
Set $f_j=\sum_{k=0}^\infty \gamma^{k} \tilde P(\gamma^{-N_j(k)})$. Conditioning on the first step $Y_1$, we get the recurrence
$$
f_j=
\begin{cases}
1+\gamma\Big[\frac{f_{j+1}+f_{j-1}}{2}\Big], & j\neq 0\\
1+\Big[\frac{f_{j+1}+f_{j-1}}{2}\Big], & j= 0.
\end{cases}
$$
The solution is $f_j=3+B(\frac{3-\sqrt{5}}{2})^j$ where $B=\frac{2}{\sqrt{5}-1}$. Since $ \sum_{j\in\Z}B(\frac{3-\sqrt{5}}{2})^j=\frac{5+\sqrt{5}}{2}$ it follows that 
$$
E[T]=3+\frac{5+\sqrt{5}}{2}q+O(q^2)
$$
as $q\to 0$. The authors thank Neal Madras for his comments about the above calculation. 

One can also approximate $E[T]$ for small $p$. For example, to first order, $E[T]$ is the mean time to reach the first $\SWE$ site on either side of the origin. So
$$
E[T]\sim\sum_{j=1}^\infty\sum_{k=1}^\infty p^2q^{j+k-1}\Big(\frac{j+k}{2}\Big)^2\sim\frac{3}{2p^2}
$$
as $p\to 0$. This can be improved on by considering additional $\SWE$ sites.

One can show in general that $v^{[2]}<0$, ie. that $E[T]<\infty$ for $q<1$. (This approach was suggested by Remco van der Hofstad).  Let $T=T(\mc{G})=\inf\{n>0:X^{[2]}_n=-1\}$, and $M_n=\#\{i<n:\mc{G}_{X_i}=\SWE\}$.  Then $v^{[2]}=-E[T]^{-1}$, and $M_n=\sum_x\ell_n(x)I_{x}$, where $I_x=I_{\{\mc{G}_x=\smallSWE\}}$ for each $x\in \Z$, and $\ell_n(x)$ denotes the local time of a {\em simple random walk} in $\Z$ (as above, the RWRE behaves as a random walk in $\Z$ killed as soon as it chooses to take a $\downarrow$ step from a $\SWE$ site).  Therefore, letting $R_n$ denote the range of simple random walk we have by Donsker-Varadhan and Markov's inequality,
\begin{align*}
P(T>n)=&E\left[\left(\frac23\right)^{M_n}\right]=E\left[\left(\frac23\right)^{\sum_x\ell_n(x)I_{x}}\right]\\
\le &P(-R_n>-n^{\frac14})+E\left[\left(\frac23\right)^{\sum_x\ell_n(x)I_{x}}\big|R_n\ge n^{1/4}\right]P(R_n\ge n^{1/4})\\
\le &e^{-Cn^{\frac13}}e^{n^\frac14}+E\left[\left(\frac23\right)^{\sum_x\ell_n(x)I_{x}}\big|R_n\ge n^{1/4}\right].
\end{align*}
The second term is at most
\begin{align*}
E\left[\left(\frac23\right)^{\sum_{x=1}^{n^{1/4}}I_{x}}\right]=E\left[\prod_{x=1}^{n^{1/4}}\left(\frac23\right)^{I_{x}}\right]=\big(\frac{2}{3}p+(1-p)\big)^{n^{1/4}}.
\end{align*}
Therefore there exists $C>0$ (independent of $p$) such that 
\begin{equation*}
P(T>n)\le e^{-Cn^{\frac13}}e^{n^\frac14}+\big(\frac{2}{3}p+(1-p)\big)^{n^{1/4}},
\end{equation*}
which is summable in $n$ for all $p>0$, whence $v^{[2]}<0$  for $p>0$.

\end{itemize}

\section{Non-monotone speeds and other calculations}

In \cite{HS_RWDRE} there are several examples of non-monotonicity, for which details of the calculations are not given. We present those here. 

\begin{itemize}
\item $\NEalpha \leftarrow$: {\it a 2-valued example in which $v^{[1]}$ is not monotone in $\alpha$}\medskip

\noindent Here $\NEalpha\vspace{2mm}$ is an environment $\gamma_1$, with $\gamma_1(e_1)=\alpha=1-\gamma_1(-e_2)$. $\gamma_2$ is the environment 
$\leftarrow$ with $\gamma_2(-e_1)=1$. We take $\mu(\{\gamma_1\})=p=1-\mu(\{\gamma_2\})$. This example appears as Example 4.6 of  \cite{HS_RWDRE}. Let $T$ be the first time we move $\uparrow$. We first calculate $E[T]$. 

Suppose that there is a $\NEalpha\vspace{2mm}$ at $(-i,0)$ for $i\ge 1$, and $\leftarrow$'s at $o$ and all points in between (a scenario with probability $p(1-p)^i$). Then $X_n$ takes $i$ steps to the left, and then oscillates between $(-i,0)$ and $(-i+1,0)$ a random number of times, before $T$ occurs. 

The other possibility is that there is a $\leftarrow$ at $(j,0)$ for $j\ge 1$, and $\NEalpha\vspace{2mm}$'s at $o$ and all points in between. This scenario has probability $(1-p)p^j$. Now $X_n$ steps right, and $T$ may occur before it reaches $(j,0)$, or it may reach $(j,0)$ and then oscillate until time $T$. The various scenarios lead to the following expression:
\begin{align*}
E[T]&=\sum_{i=1}^\infty p(1-p)^i\sum_{k=0}^\infty \alpha^k(1-\alpha)[i+2k+1]\\
&\qquad\qquad+\sum_{j=1}^\infty (1-p)p^j\Big(\sum_{k=0}^{j-2}\alpha^k(1-\alpha)[k+1] + \sum_{k=0}^\infty \alpha^{j+k-1}(1-\alpha)[j+2k]\Big)\\
&=\sum_{i=1}^\infty p(1-p)^i(i+\frac{1+\alpha}{1-\alpha})+\sum_{j=1}^\infty (1-p)p^j\Big((1-\alpha)\frac{d}{d\alpha}\frac{1-\alpha^j}{1-\alpha}+j\alpha^{j-1}+\frac{2\alpha^j}{1-\alpha}\Big)\\
&=\frac{1+\alpha}{1-\alpha}(1-p)+\frac{1-p}{p}+\sum_{j=1}^\infty (1-p)p^j\frac{1+\alpha^j}{1-\alpha}\\
&=\frac{1+\alpha}{1-\alpha}(1-p)+\frac{1-p}{p}+\frac{p}{1-\alpha}+\frac{(1-p)p\alpha}{(1-\alpha)(1-p\alpha)}\\
&=\frac{\alpha}{1-\alpha}+p+\frac{1-p}{p}+\frac{(1-p)}{(1-\alpha)(1-p\alpha)}.
\end{align*}
We expect this to be decreasing in $p$ and increasing in $\alpha$. Rewriting the first term as $\frac{1}{1-\alpha}-1$ we see that indeed, each term is increasing in $\alpha$. It can also be rewritten as
$$
\frac{\alpha}{1-\alpha}+\frac{1}{\alpha(1-\alpha)}+p+\frac{1}{p}-1-\frac{1}{1-p\alpha}(\frac{1}{\alpha}-1)
$$
and both $p+\frac{1}{p}$ and $-1/(1-p\alpha)$ are indeed decreasing in $p$, so our expectations are realized.
This gives us that $v^{[2]}=1/E[T]$ is monotone in both parameters. 

The martingale argument given above shows that 
$$
v^{[j]}E[T]=\mc{L}_1f_jN_1 + \mc{L}_2f_jN_2
$$
where the $\mc{L}_i$ are the generators for the two environments, $f_j(x)=x^{[j]}$, and $N_i$ is the expected number of visits to environment $i$ before $T$. In particular, $\mc{L}_1f_1=\alpha$, $\mc{L}_1f_2=1-\alpha$, $\mc{L}_2f_1=-1$, $\mc{L}_2f_2=0$, so 
$$
v^{[1]}E[T]=\alpha N_1-N_2, \qquad 1=v^{[2]}E[T]=(1-\alpha) N_1, \qquad N_1+N_2=E[T].
$$
Therefore $N_1=\frac{1}{1-\alpha}$ which implies that
\begin{align*}
v^{[1]}&=\frac{1+\alpha}{1-\alpha}\frac{1}{E[T]}-1\\
&=(1+\alpha)\left(\alpha+(1-\alpha)\Big(p+\frac1p-1\Big)+\frac{(1-p)}{(1-p\alpha)}\right)^{-1}-1.
\end{align*}

When $\alpha=1$, $v^{[1]}=0$. When $\alpha=0$, $v^{[1]}=p-1<0$. Graphing $v^{[1]}$, we see that it is increasing in $\alpha$ for $p\le 0.5$. But for $p\ge .58$,  $v^{[1]}$ has a positive interior maximum.

Likewise, consider the fraction of sites $\pi^{\sss 1}$ which are of type $\NEalpha$, namely $N_1/E[T]=1/W$, where 
$$
W=(1-\alpha)E[T] = p+\frac1p-1+\alpha(2-(p+\frac1p))+\frac{1-p}{1-\alpha p}.
$$
When plotted, this function is not monotone in $\alpha$, at least when $p\approx\frac12$. 

When $\alpha=0$ we have $W=\frac{1}{p}$. When $\alpha=1$ we have $W=2$, and clearly $p=\frac12$ is where these two values agree. In other words, for small $p$ the fraction of $\NEalpha\vspace{2mm}$ is increasing with $\alpha$, whereas for large $p$ it is decreasing with $\alpha$.  In other words, for small $p$ an increase in $\alpha$ means we tend to linger more around the $\NEalpha\vspace{2mm}$ we've managed to find. Whereas when $p$ is big, increasing $\alpha$ means we get stuck for a longer time, which brings the fraction of $\NEalpha$ we see down towards $\frac12$. \medskip

\item $\NEalpha\vspace{2mm} \NEbeta\leftarrow:$
{\it a 3-valued example in which the frequency $v^{[1]}$ is not monotone in $p_1$.}

This example contains 4 parameters $(\alpha,\beta,p,q)\in [0,1]^4$.
We let $\gamma^1(e_1)=\alpha=1-\gamma^1(e_2)$, $\gamma^2(e_1)=\beta=1-\gamma^2(e_2)$, and $\gamma^3(-e_1)=1$.  
We set $\mu(\{\gamma^3\})=p$, $\mu(\{\gamma^1\})=(1-p)q$, and $\mu(\{\gamma^2\})=(1-p)(1-q)$. The claim is that there exist $\alpha,\beta,p\in (0,1)$ such that $v^{[1]}=v^{[1]}(\alpha,\beta,p,q)$ is not monotone increasing in $q$. This appears as Example 4.4 of \cite{HS_RWDRE}

We need to calculate $E[T]$ and $E[X^{[1]}_T]$.  

Let $Y^-$ be the number of consecutive $\leftarrow$ starting from the origin (going to the left).  Then for $k\ge 0$, $\mP(Y^-=k)=p^{k}(1-p)$.  Let $Y^+$ be the number of consecutive $\NE$ (of either type) starting from the origin (going to the right).  Then for $k\ge 0$, $\mP(Y^+=k)=(1-p)^{k}p$.  Let $\eta=q\alpha+(1-q)\beta$ and $\xi=\frac{q\alpha}{1-\alpha}+\frac{(1-q)\beta}{1-\beta}$.  Let $\Gamma$ be a random variable that is equal to $1-\alpha$ with probability $q$ and $1-\beta$ otherwise.  Let $G\ge 0$ be a random variable that (conditional on $\Gamma$) is Geometric with parameter $\Gamma$.

Let $\Delta=X_T^{[1]}$.  Then 
\[\mE[\Delta]=\mE[\Delta\mathbbm{1}_{Y^->0}]+\mE[\Delta\mathbbm{1}_{Y^+>0}],\]
where the first term on the right is 
\[-\mE[Y^-\mathbbm{1}_{Y^->0}]=-\mE[Y^-]=\frac{-p}{1-p}.\]
Similarly
\[\mE[T]=\mE[T\mathbbm{1}_{Y^->0}]+\mE[T\mathbbm{1}_{Y^+>0}],\]
where the first term on the right is equal to
\[\mE[(1+Y^-+2G)\mathbbm{1}_{Y^->0}]=p(1+2\xi)+\mE[Y^-\mathbbm{1}_{Y^->0}]=p(1+2\xi)+\frac{p}{1-p}.\]

Consider now $\mE[\Delta\mathbbm{1}_{Y^+>0}]$ and $\mE[T\mathbbm{1}_{Y^+>0}]$.  If $Y^+=k>0$, the walker has $k-1$ independent opportunities to go up, each with probability $\eta$, before reaching the last $\NE$ site before the first $\leftarrow$ site.  Under this conditioning, the probability that $\Delta=j\in [0,k-2]$ (and $T=j+1$) is then $\eta^{j}(1-\eta)$.  The probability that $T=k-1+2m+1=k+2m$ is 
\[\eta^{k-1}\left[q\alpha^m(1-\alpha)           +(1-q)\beta^m(1-\beta)   \right].\]
The probability that $\Delta=k-1$ is $\eta^{k-1}$ (note that this is also the sum over $m\ge 0$ of the above).  Thus, using $\sum_{k=1}^{\infty}kr^{k-1}=(1-r)^{-2}$ and letting $a=(1-p)\eta$ we have

\begin{align}
\mE[\Delta\mathbbm{1}_{Y^+>0}]=&\sum_{k=1}^{\infty}(1-p)^{k}p\left((k-1)\eta^{k-1} +\sum_{j=0}^{k-2}j\eta^{j}(1-\eta)\right)\\
=&\frac{p(1-p)^2\eta}{(1-a)^2}+\frac{\eta (1-\eta)(1-p)^3}{(1-a)^2}=\frac{(1-p)^2\eta}{1-a}
\end{align}
Moreover,
\begin{align}
\mE[T\mathbbm{1}_{Y^+>0}]=&\sum_{k=1}^{\infty}(1-p)^{k}p\Bigg(\sum_{m=0}^{\infty}(k+2m)\eta^{k-1}\left[q\alpha^m(1-\alpha)           +(1-q)\beta^m(1-\beta)   \right]\\ &+\sum_{j=0}^{k-2}(j+1)\eta^{j}(1-\eta)\Bigg)\\
=&(1-p)p\left[q(1-\alpha)h(\alpha)+(1-q)(1-\beta)  h(\beta) \right]  +\frac{(1-p)^2(1-\eta)}{(1-a)^2},
\end{align}
where 
\[h(x)=\sum_{k=1}^{\infty}a^{k-1}\sum_{m=0}^{\infty}(k+2m)x^m=\frac{1}{(1-a)^2(1-x)}+\frac{2x}{(1-a)(1-x)^2}.  \]

Thus,
\begin{align}
\mE[T\mathbbm{1}_{Y^+>0}]=\frac{1-p}{1-a}\left[1+2p\xi\right].
\end{align}


Finally we have that
\begin{align}
\mE[\Delta]=&\frac{(1-p)^2\eta}{1-a}-\frac{p}{1-p}\\
\mE[T]=&p(1+2\xi)+\frac{p}{1-p}    +\frac{1-p}{1-a}\left[1+2p\xi\right],
\end{align}
where $a=(1-p)\eta=(1-p)(q\alpha+(1-q)\beta)$ and $\xi=\frac{q\alpha}{1-\alpha}+\frac{(1-q)\beta}{1-\beta}$.

\end{itemize}

\bibliographystyle{plain}

\end{document}